\documentclass[a4paper, 11pt]{amsart}

\usepackage[mathscr]{eucal}
\usepackage{amsmath}
\usepackage{amssymb}
\usepackage{amsfonts}
\usepackage{amsthm}

\usepackage{cases}

\usepackage[mathscr]{eucal}
\usepackage{amsmath}
\usepackage{amssymb}
\usepackage{amsfonts}
\usepackage{amsthm}
\usepackage[dvips]{graphics, color}
\theoremstyle{plain}
\newtheorem{theorem}{Theorem}[section]

\newtheorem{corollary}{Corollary}[section]

\newtheorem*{problem*}{Problem}

\theoremstyle{definition}
\newtheorem{example}{Example}[section]
\newtheorem{remark}{Remark}[section]
\newtheorem*{remark*}{Remark}

\theoremstyle{remark}

\numberwithin{equation}{section}

\def\<{\left<} \def\>{\right>}

\def\proof{\noindent{\it Proof. }}
\def\bea{\begin{eqnarray} }
\def\eea{\end{eqnarray} }
\def\be{\begin{equation} }
\def\ee{\end{equation} }
\def\qed{\ifhmode\unskip\nobreak\fi\ifmmode\ifinner\else\hskip5pt
\fi\fi\hbox{\hskip5 pt \vrule width4 pt height6 pt depth1.5 pt \hskip1pt }}

\begin{document}
\title[]{Real hypersurfaces in the complex projective plane satisfying an equality involving $\delta(2)$}
\author[]{Toru Sasahara}
\address{
Center for Liberal Arts and Sciences, 
Hachinohe Institute of Technology, 
Hachinohe, Aomori 031-8501, Japan}
\email{sasahara@hi-tech.ac.jp}


\begin{abstract} 
It was proved in Chen's paper \cite{chen} that
  every real hypersurface in the complex projective plane of constant holomorphic sectional curvature $4$ satisfies
$$ \delta(2)\leq \frac{9}{4}H^2+5,$$
where $H$ is the mean curvature and $\delta(2)$ is a $\delta$-invariant introduced by him.
In this paper, we  
study non-Hopf real hypersurfaces satisfying the equality case of the 
inequality under the condition that
 the mean curvature is constant along each integral curve of the Reeb vector field.
 We describe how to obtain all such hypersurfaces.
\end{abstract}

\keywords{real hypersurfaces, ruled, $\delta(2)$-ideal,
 complex projective plane.}

\subjclass[2010]{Primary 53C42; Secondary 53B25.}

\maketitle

 \section{Introduction}


For a  Riemannian $m$-manifold $M$ with $m>2$, Chen \cite{chen2} introduced in the early 1990s the following
invariant:
$$\delta(2)(p)=\tau(p)-\inf\{K(\pi)\ |\ \pi\  \text{is a plane in}\ T_pM\},$$
where  $\tau$ is the scalar curvature and  $K(\pi)$ is the sectional curvature of $\pi$.
If $m=3$, then $\delta(2)(p)$ is equal to  the maximum Ricci curvature function $\overline{Ric}$
on $M$ defined by
$\overline{Ric}(p)=\max\{S(X, X)\ |\ X\in T_pM, \ ||X||=1\}$,
where $S$ is the Ricci tensor.  For general $\delta$-invariants,
see \cite{chen6} for details. 

It was proved in \cite{chen} that every real hypersurface in the complex projective space
$\mathbb{C}P^n$ of complex dimension $n$
and constant holomorphic sectional curvature $4$ satisfies
\be
 \delta(2)\leq \frac{(2n-1)^2(2n-3)}{4(n-1)}H^2+2n^2-3,\label{ideal}
\ee
where $H$ denotes the mean curvature.
A real hypersurface in $\mathbb{C}P^n$ is said to be $\delta(2)$-{\it ideal} if it attains equality in (\ref{ideal}) at each point.
 Chen \cite{chen} completely classified $\delta(2)$-ideal Hopf real hypersurfaces in $\mathbb{C}P^n$. 
In \cite{sa},  the author proved  that a non-Hopf
real hypersurface   with constant mean curvature in $\mathbb{C}P^2$ is $\delta(2)$-ideal 
 if and only if it is a minimal ruled real hypersurface.
In this paper, we  
 classify $\delta(2)$-ideal non-Hopf real hypersurfaces in  $\mathbb{C}P^2$ 
 whose mean curvature is constant along each integral curve of the Reeb vector field.


 
  
%
%

\section{Preliminaries}
Let $M$ be a   real hypersurface in the complex projective space  $\mathbb{C}P^n$.
 We denote by $J$ the almost 
complex structure of $\mathbb{C}P^n$.  
For a   unit normal vector field $N$, the vector field on $M$ defined by  $\xi=-JN$ is called the {\it Reeb
	vector field}.
If $\xi$ is a principal curvature vector at every point of  $M$, then $M$ is said to be {\it Hopf}.

Let $\mathcal{H}$ be the holomorphic distribution defined by $\mathcal{H}=\bigcup_{p\in M}\{X\in T_pM\ |\ \<X, \xi\>=0\}$, where $\<\cdot,\cdot \>$ denotes the metric of $\mathbb{C}P^n$.
If $\mathcal{H}$ is integrable and each leaf of its maximal integral manifolds is 
a totally geodesic complex hypersurface,
then $M$ is said to be  {\it ruled}.

Denote by $\nabla$ and
 $\tilde\nabla$ the Levi-Civita connections on $M$ and $\mathbb{C}P^n$, respectively. 
  The
 Gauss and Weingarten formulas are respectively given by
\be
 \begin{split}
 \tilde \nabla_XY&= \nabla_XY+\<AX, Y\>N, \label{gawe}\\
 \tilde\nabla_X N&= -AX
 \end{split}\nonumber
\ee
 for tangent vector fields $X$, $Y$ and a unit normal vector field $N$,
 where $A$ is the shape operator with respect to $N$.
 The function $H={\rm tr}A/(2n-1)$ is called the  {\it mean curvature}.
 If it vanishes identically, then $M$ is said to be  {\it minimal}.
 
For any  vector field  $X$ tangent to $M$,  we denote the tangential component of $JX$ by $\phi X$.
Then by the Gauss and  Weingarten formulas, we have
\be
\nabla_{X}\xi=\phi AX. \label{PA}
\ee
 We denote by $R$ the Riemannian curvature tensor of $M$. Then,
 the equations of Gauss  and Codazzi are respectively given by
 \begin{align}
 &R(X, Y)Z=\<Y, Z\>X-\<X, Z\>Y+\<\phi Y, Z\>\phi X
 -\<\phi X, Z\>\phi Y  \label{ga}\\
 &\hskip60pt -2\<\phi X, Y\>\phi Z
   +\<AY, Z\>AX-\<AX, Z\>AY,\nonumber\\
  &(\nabla_XA)Y-(\nabla_YA)X=\<X, \xi\>\phi Y-\<Y, \xi\>\phi X-2\<\phi X, Y\>\xi.\label{co}
 \end{align}       
 
 \section{$\delta(2)$-ideal real hypersurfaces} 
Applying  \cite[Theorem 5]{chen} to real hypersurfaces in $\mathbb{C}P^n$, we have the following general inequality.
\begin{theorem}\label{thm1}
	Let $M$ be a real hypersurface in 	$\mathbb{C}P^n$.
For any point $p\in M$ and any plane $\pi\subset T_pM$, we  have
\be
\tau-K(\pi)\leq\frac{(2n-1)^2(2n-3)}{4(n-1)}H^2+2n^2-3-3\<Je_1, e_2\>^2 \label{id},
\ee
 where $\{e_1, e_2\}$ is an orthonormal basis of $\pi$.	
The equality sign in $(\ref{id})$ holds at a point $p\in M$ 
if and only if there exists an orthonormal basis $\{e_1, e_2,  \ldots, e_{2n-1}\}$ at $p$ such that 
the shape operator  at $p$ is represented by a matrix 
\be
A= \left(
\begin{array}{ccccc}
	\alpha & \beta & 0 & \ldots & 0 \\
	\beta & \gamma & 0 & \ldots & 0 \\ 
	0 & 0 & \mu &\dots & 0 \\
	\vdots & \vdots & \vdots & \ddots & \vdots \\
	0 & 0 & 0 & \ldots & \mu
\end{array}
\right), \label{A} 
\ee
where $\alpha+\gamma=\mu$.
%
	\end{theorem}

The following Corollary immediately follows from Theorem \ref{thm1}.
\begin{corollary}[\cite{chen}]\label{cor1}
	Let $M$ be a real  hypersurface 
	in $\mathbb{C}P^n$. Then, we have
	\be
	\delta(2)\leq \frac{(2n-1)^2(2n-3)}{4(n-1)}H^2+2n^2-3 \label{ideal2}
	\ee
	at each point of $M$.
	The equality sign in $(\ref{ideal2})$ holds at a point $p\in M$ 
	if and only if there exists an orthonormal basis $\{e_1, e_2,  \ldots, e_{2n-1}\}$ at $p$ such that
	
	${\rm (1)}$ $\<Je_1, e_2\>=0$,  
	
	${\rm (2)}$ $K(e_1\wedge e_2)={\rm inf} K$,
	
	${\rm (3)}$
	the shape operator  at $p$ is represented by a matrix {\rm (\ref{A})} with $\alpha+\gamma=\mu$.
\end{corollary}

\begin{remark}\label{rm1}
It follows from (\ref{id}) that	if (1) and (3) in Corollary \ref{cor1} hold, then (2) is automatically satisfied. 
	\end{remark}

A real hypersurface in $\mathbb{C}P^n$ is said to be $\delta(2)$-{\it ideal} if it attains equality in (\ref{ideal2}) at each point.
In \cite{chen}, Chen 
proved that a Hopf real hypersurface  in  $\mathbb{C}P^n$ is $\delta(2)$-ideal
if and only if it is an open part of one of the following hypersurfaces:
(i) a geodesic sphere with radius $\pi/4$ in  $\mathbb{C}P^n$,
(ii)  a tubular hypersurface with radius 
$r=\tan^{-1}((1+\sqrt{5}-\sqrt{2+2\sqrt{5}})/2)$ over a
complex quadric curve $Q_1$  in  $\mathbb{C}P^2$.

We now present  a class of $\delta(2)$-ideal non-Hopf hypersurfaces in $\mathbb{C}P^2$.
\begin{example}\label{ex1}
	Suppose that $\alpha(s)$, $\beta(s)$, $\gamma(s)$ and $\mu(s)$ satisfy
	\bea
	&&\begin{split}\label{2-hopf}
		&\alpha^{\prime}=\beta(\alpha+\gamma-3\mu),\\
		& \beta^{\prime}=\beta^2+\gamma^2+\mu(\alpha-2\gamma)+1,\\
		&\gamma^{\prime}=\frac{(\gamma-\mu)(\gamma^2-\alpha\gamma-1)}{\beta}
		+\beta(2\gamma+\mu),
	\end{split}
	\eea
	on an open interval $I\subset\mathbb{R}$, 
	where $\beta(s)$ are nowhere zero.
	According to Theorem 5 in \cite{ivey},
	there exists a smooth immersion $\Phi: I\times\mathbb{R}^2\rightarrow\mathbb{C}P^2$
	determining a non-Hopf real hypersurface  in $\mathbb{C}P^2$,
	such that
	the shape operator $A$ is 
	represented by {\rm (\ref{A})}
	with respect to an orthonormal frame field $\{\xi, X, \phi X\}$, where $\phi X=\partial/\partial s$. 
	The distribution  $\mathcal{D}$ spanned by $\xi$ and $X$ is integrable, and $\Phi$ maps the $\mathbb{R}^2$-factors onto the $\mathcal{D}$-leaves. 
	Clearly, the mean curvature of the hypersurface is constant along each integral curve of the Reeb vector field.
	
	 If $\alpha+\gamma=\mu$ on $I$, then Corollary \ref{cor1} and Remark \ref{rm1} imply
	that $\Phi$ is $\delta(2)$-ideal. In particular, if $\alpha=\gamma=\mu=0$ on $I$, then ${\rm tr}A=0$ and $\<AX, Y\>=0$ for
	any tangent vector field $X$, $Y$ on $M$ orthogonal to $\xi$, and hence $\Phi$ is minimal ruled 
	(see \cite[p.445]{ce} and \cite{kim1}).
	\end{example}

\begin{remark}
	Substitution of $\alpha+\gamma=\mu$ into (\ref{2-hopf}) gives a autonomous system.
	It follows from Picard's theorem 
	that for given initial values $\alpha(s_0)=\alpha_0$, $\beta(s_0)=\beta_0$, $\gamma(s_0)=\gamma_0$ with $\beta_0\ne 0$ and $\alpha_0+\gamma_0\ne 0$, 
	the initial value problem of   (\ref{2-hopf}) with $\alpha+\gamma=\mu$
	has a unique solution satisfying $\beta\ne 0$ and $\alpha+\gamma\ne 0$
	on some open interval containing $s_0$.
	Therefore, there exist infinity many $\delta(2)$-ideal real hypersurfaces in 
	$\mathbb{C}P^2$ which are non-Hopf and non-minimal.
\end{remark}

\begin{remark}
	Let $M$ be a real hypersurface in the complex hyperbolic space
	$\mathbb{C}H^n$ of constant
	holomorphic sectional curvature $-4$. Then we have
	\be
	\delta(2) \leq \frac{(2n-1)^2(2n-3)}{4(n-1)}H^2+6-2n^2. \nonumber
	\ee
	The equality sign of the inequality holds identically if and only if
	$M$ is an open part of the horosphere in $\mathbb{C}H^2$ (see \cite{chen}). 
\end{remark}

\section{Main result}
The following theorem is the main result of this paper.
\begin{theorem}\label{Non-Hopf}
	Let $M$ be a  $\delta(2)$-ideal non-Hopf real hypersurface  in  $\mathbb{C}P^2$.
	If the mean curvature is constant along each integral curve of the Reeb vector field, then
	$M$  is locally obtained by the construction described in Example $\ref{ex1}$.
\end{theorem}

\proof
  Let $M$ be a $\delta(2)$-ideal non-Hopf real hypersurface  in $\mathbb{C}P^2$. 
Let $\{e_1, e_2, e_3\}$ be a local orthonormal frame field described in Corollary \ref{cor1}. 
We put $\xi=pe_1+qe_2+re_3$ for
some functions $p$, $q$ and $r$.
It follows from  $\<Je_1, e_2\>=0$ that $r\<Je_3, e_1\>=r\<Je_3, e_2\>=0$.
If $r\ne 0$, then $\xi=e_3$. However, this contradicts  $\<Je_1, e_2\>=0$.
Hence,  $r=0$ holds, that is, $\xi$ lies in Span$\{e_1, e_2\}$.
We may assume that $e_1=\xi$ and $Je_2=e_3$. 
 From  (3) of Corollary \ref{cor1}, we see that  the shape operator satisfies the following:
\begin{align} 
& A\xi=(\mu-\gamma)\xi+\beta e_2,\ \  Ae_2=\gamma e_2+\beta\xi, \ \ Ae_3
=\mu e_3.\label{A2}
 \end{align}

Let $\Omega$ be an open set where $\beta\ne 0$. We work in $\Omega$.
  Using (\ref{PA}) and (\ref{A2}), we get
\be
\nabla_{e_2}\xi=\gamma e_3, \ \ \nabla_{e_3}\xi=-\mu e_2, \ \ \nabla_{\xi}\xi=\beta e_3. \label{nab1}
\ee 
Since $\<\nabla e_i, e_j\>=-\<\nabla e_j, e_i\>$ holds, by (\ref{nab1}) we have
  \be
\begin{split}
&   \nabla_{e_2}e_2 =\kappa_1e_3, \ \ \nabla_{e_3}e_2=\kappa_2e_3+\mu\xi, 
\ \ \nabla_{\xi}e_2 =\kappa_3e_3, \\
& \nabla_{e_2}e_3 =-\kappa_1e_2-\gamma\xi, \ \ \nabla_{e_3}e_3=-\kappa_2e_2, 
\ \ \nabla_{\xi}e_3=
-\kappa_3e_2-\beta\xi \label{nab2}
\end{split}
\ee
for some functions  $\kappa_1$, $\kappa_2$ and $\kappa_3$.

 
 Assume that the mean curvature $H=\mu/3$ is constant along each integral curve of the Reeb vector field $\xi$,
 that is, 
 \be
 \xi\mu=0. \label{assum}
 \ee
From (\ref{A2}), (\ref{nab1}), (\ref{nab2}) and    the equation  (\ref{co}) of Codazzi, it follows  that
\begin{align}
e_2\mu &=0, \label{cd1}\\
e_3\gamma &=(\gamma-\mu)\kappa_1+\beta(\gamma+2\mu),\label{cd2}\\
e_3\beta &=-\gamma^2+\beta\kappa_1-2\gamma\mu+\mu^2+2\label{cd3},\\
e_2\beta &=\xi\gamma,\label{cd4} \\
e_2\gamma &=-\xi\beta, \label{cd5} \\
\beta\kappa_1+(\mu-\gamma)\kappa_3 &=\beta^2+\gamma^2-1,\label{cd6}\\
\kappa_2 &=0, \label{cd7}\\
e_3(\mu-\gamma) &=\beta(\kappa_3-2\mu-\gamma). \label{cd8} 
\end{align}
Taking into account (\ref{cd7}), the equation (\ref{ga}) of Gauss for $\<R(e_2, e_3)e_3, e_2\>$ and $\<R(\xi, e_2)e_3, e_2\>$ yields
\begin{align}
e_3\kappa_1&=2\mu\gamma+\kappa_1^2+(\gamma+\mu)\kappa_3+4, \label{ga1}\\
\xi\kappa_1 &=e_2\kappa_3.\label{ga2} 
\end{align}
Using (\ref{nab1}), (\ref{nab2}),  (\ref{assum}) and (\ref{cd1}) we have
\be
0=[e_2, \xi]\mu=(\nabla_{e_2}\xi-\nabla_{\xi}e_2)\mu=(\gamma-\kappa_3)e_3\mu.\label{eq00}
\ee
Thus, we obtain that   $\gamma=\kappa_3$ or $e_3\mu=0$.

{\bf Case (a)}: $e_3\mu=0$ on an open subset ${U}\subset \Omega$.
 In this case, combining (\ref{assum}) and  (\ref{cd1})  implies that 
$\mu$ is constant, that is, the mean curvature is constant on ${U}$. Hence, by virtue of  \cite[Theorem 1.2]{sa}, 
we conclude that 
$U$ is minimal ruled.

{\bf Case (b)}: $\gamma=\kappa_3$ on an open subset ${V}\subset \Omega$. 
In this case,  
since $\nabla_{e_2}\xi-\nabla_{\xi}e_2=0$ holds, the distribution
  $\mathcal{D}$ spanned by $\xi$ and $e_2$ is integrable.
Eliminating $e_3\gamma$ from (\ref{cd2}) and (\ref{cd8}), we obtain
\be 
e_3\mu=(\gamma-\mu)\kappa_1+\beta\gamma. \label{eq2.1}
\ee
Equations (\ref{cd6}) and (\ref{ga1}) become
\begin{align}
\beta\kappa_1 &=\beta^2+2\gamma^2-\mu\gamma-1, \label{eq2}\\
e_3\kappa_1 &=\kappa_1^2+\gamma^2+3\gamma\mu+4, \label{eq2.2}
\end{align}
respectively. From (\ref{cd5}) and  (\ref{ga2}), it follows that
\be
\xi\kappa_1=-\xi\beta.\label{eq1}
\ee
 Elimination of  $\kappa_1$ from (\ref{cd3}) and (\ref{eq2}) leads to
\be
e_3\beta=\beta^2+\gamma^2-3\gamma\mu+\mu^2+1. \label{eq3}
\ee
Using (\ref{nab1}), (\ref{nab2}), (\ref{cd2}), (\ref{cd4}), (\ref{cd7}), (\ref{eq1}) and  (\ref{eq3}), we have
 the following:
\begin{align}
e_3(\xi\beta) &=(\nabla_{e_3}\xi-\nabla_{\xi}e_3)\beta+\xi(e_3\beta) \nonumber \\
&=(\gamma-\mu)\xi\gamma+\beta(\xi\beta)+\xi(\beta^2+\gamma^2-3\gamma\mu+\mu^2+1)\nonumber\\
&=3\beta(\xi\beta)+(3\gamma-4\mu)\xi\gamma, \label{eq3.1}\\
e_3(\xi\gamma) &=(\nabla_{e_3}\xi-\nabla_{\xi}e_3)\gamma+\xi(e_3\gamma)  \nonumber \\
&=(\mu-\gamma)\xi\beta+\beta(\xi\gamma)+\xi[(\gamma-\mu)\kappa_1+\beta(\gamma+2\mu)]\nonumber\\
&=(4\mu-\gamma)\xi\beta+(2\beta+\kappa_1)\xi\gamma.\label{eq3.2}
\end{align}

Differentiating (\ref{eq2}) with respect to 
$\xi$, and using (\ref{assum}) and (\ref{eq1}),  we obtain
\be
(\kappa_1-3\beta)\xi\beta+(\mu-4\gamma)\xi\gamma=0. \label{eq4}
\ee
Moreover, differentiating (\ref{eq4}) with respect to  $e_3$, 
we have
\be 
(e_3\kappa_1-3e_3\beta)\xi\beta+(\kappa_1-3\beta)e_3(\xi\beta)+(e_3\mu-4e_3\gamma)\xi\gamma
+(\mu-4\gamma)e_3(\xi\gamma)=0.\label{eq4.1}
\ee
Substitution of  (\ref{cd2}), (\ref{eq2.1}), (\ref{eq2.2}),  (\ref{eq3}), (\ref{eq3.1}) and (\ref{eq3.2}) 
into (\ref{eq4.1}) gives
\be
(\kappa_1^2-12\beta^2+2\gamma^2+\mu^2-5\mu\gamma+3\beta\kappa_1+1)\xi\beta+
(6\beta\mu-20\beta\gamma-4\gamma\kappa_1)\xi\gamma=0. \label{eq5}
\ee
Equations (\ref{eq4}) and (\ref{eq5}) could be rewritten as 
\be
    \begin{pmatrix}
       a_{11}&a_{12} \\
       a_{21}& a_{22}\\
   \end{pmatrix}
  \begin{pmatrix}
  \xi\beta \\
  \xi\gamma
  \end{pmatrix}=\begin{pmatrix}
  0 \\
  0
  \end{pmatrix}, \label{eq6}
\ee
where the components of the  square matrix are given by
\begin{align*}
& a_{11}=\kappa_1-3\beta, \\
&a_{12}=\mu-4\gamma,\\
&a_{21}=\kappa_1^2-12\beta^2+2\gamma^2+\mu^2-5\mu\gamma+3\beta\kappa_1+1, \\
&a_{22}=6\beta\mu-20\beta\gamma-4\gamma\kappa_1.
\end{align*}
We divide Case (b) into two subcases.

{\bf Case (b.1)}: $a_{11}a_{22}-a_{21}a_{12}\ne 0$ on an open neighborhood ${V}_1$
of a point $p\in{V}$.
In this case, by (\ref{eq6}), we have $\xi\beta=\xi\gamma=0$. It follows from (\ref{cd4}) and (\ref{cd5}) that
$e_2\beta=e_2\gamma=0$. This, together with  (\ref{assum}) and (\ref{cd1}), 
 implies that 
all the components of the shape operator $A$ are constant along the 
$\mathcal{D}$-leaves.  Moreover, equations (\ref{cd2}), (\ref{cd3}) and (\ref{cd8})
imply that  (\ref{2-hopf}) with $\alpha+\gamma=\mu$, where $d/ds$ stands for the derivative with respect to $e_3$.
Note that the existence of such a hypersurface is guaranteed by Example \ref{ex1}.

{\bf Case (b.2)}:  $a_{11}a_{22}-a_{21}a_{12}=0$ on an open neighborhood ${V}_2$
of a point $p\in{V}$.
In this case, eliminating $\kappa_1$ from this condition and (\ref{eq2}) yields
\be
p_1(\gamma, \mu)\omega^2+p_2(\gamma, \mu)\omega+p_3(\gamma, \mu)=0,\label{eq7}
\ee
where $\omega=\beta^2$, and $p_i$ are polynomials  given by
\begin{align*}
&p_1=16\gamma-4\mu, \\
&p_2=16\gamma^3-24\gamma^2\mu+8\gamma\mu^2-\mu^3-2\mu,\\
&p_3=-\mu(2\gamma^2-\gamma\mu-1)^2.
\end{align*}
Differentiating (\ref{eq7}) with respect to $e_3$, and  using
(\ref{cd2}), (\ref{eq2.1})  and 
(\ref{eq3}), we obtain
\be
\begin{split}
&\kappa_1(12\beta^4-12\beta^4\mu+24\beta^2\gamma^3
-56\beta^2\gamma^2\mu+37\beta^2\gamma\mu^2-2\beta^2\gamma \\
&-5\beta^2\mu^3+2\beta^2\mu-4\gamma^5-4\gamma^4\mu+17\gamma^3\mu^2+4\gamma^3 \\
&-11\gamma^2\mu^3+2\gamma\mu^4-6\gamma\mu^2-\gamma+2\mu^3+\mu) \\
&+76\beta^5\gamma+16\beta^5\mu+120\beta^3\gamma^3-
192\beta^3\gamma^2\mu+37\beta^3\gamma\mu^2 \\
&+62\beta^3\gamma-2\beta^3\mu^3-20\beta^3\mu+28\beta\gamma^5
-152\beta\gamma^4\mu\\
&+169\beta\gamma^3\mu^2+36\beta\gamma^3
-76\beta\gamma^2\mu^3-48\beta\gamma^2\mu+18\beta\gamma\mu^4 \\
&+42\beta\gamma\mu^2-\beta\gamma-2\beta\mu^5-10\beta\mu^3-4\beta\mu
\end{split}\label{kappa1}
\ee
Eliminating $\kappa_1$ from (\ref{kappa1}) and 
(\ref{eq2}), we get
\be
q_1(\gamma, \mu)\omega^3+q_2(\gamma, \mu)\omega^2+q_3(\gamma, \mu)\omega
+q_4(\gamma, \mu)=0, \label{eq8}
\ee
 where $\omega=\beta^2$, and $q_i$ are polynomials given by
 \begin{align*}
  q_1=& 88\gamma+4\mu,\\
 q_2=& 168\gamma^3-284\gamma^2\mu+86\gamma\mu^2+48\gamma-7\mu^3-6\mu,\\
 q_3=& 72\gamma^5-292\gamma^4\mu+316\gamma^3\mu^2+12\gamma^3-134\gamma^2\mu^3\\ 
 &+14\gamma^2\mu+25\gamma\mu^4-3\gamma\mu^2-2\mu^5-3\mu^3-5\mu,\\
 q_4=& (\mu-\gamma)(2\gamma^2-\gamma\mu-1)^2(2\gamma^2+5\gamma\mu-2\mu^2-1).
 \end{align*}
The resultant $R_1(\gamma, \mu)$ of the left-hand sides of (\ref{eq7}) and (\ref{eq8}) with respect to $\omega$
is found to be the following polynomial:
\be 
R_1(\gamma, \mu)=32(4\gamma-\mu)(2\gamma^2-\gamma\mu-1)^3
\Bigl(1536\gamma^8+\sum_{i=0}^{7}g_{i}(\mu)\gamma^i\Bigr),\nonumber
\ee
where $g_i$ are polynomials given by
\begin{align*}
& g_0=3\mu^8-8\mu^6+6\mu^4,\\
& g_1=-78\mu^7+60\mu^5+52\mu^3,\\
& g_2=720\mu^6+204\mu^4+160\mu^2,\\
& g_3=-3040\mu^5-1016\mu^3+112\mu,\\
& g_4=6752\mu^4+576\mu^2+32,\\
& g_5=-9152\mu^3-608\mu,\\
& g_6=8256\mu^2-192,\\
& g_7=-4480\mu.
\end{align*}

{\bf Case (b.2.i)}: $4\gamma-\mu=0$ on an open subset $V_{21}\subset V_2$. 
 Differentiating this condition with respect to $e_3$,
 and using (\ref{cd2}),    (\ref{eq2.1}) and (\ref{eq2}), we obtain
 \be 
 6\gamma^3-9\gamma^2\mu+3(\mu^2+2\beta^2-1)\gamma+(5\beta^2+3)\mu=0.\nonumber
 \ee
Eliminating $\gamma$ from this equation and $4\gamma-\mu=0$ yields
$$\mu(9\mu^2+208\beta^2+72)=0,$$
 which shows that $\mu=\gamma=0$ and hence $V_{21}$ is minimal
ruled.

{\bf Case (b.2.ii)}: $2\gamma^2-\gamma\mu-1=0$ on an open subset $V_{22}\subset V_2$.
Differentiating this condition with respect to $e_3$,
and using (\ref{cd2}),    (\ref{eq2.1}) and (\ref{eq2}), 
we get
$$6\gamma^4-11\gamma^3\mu+(6\mu^2+6\beta^2-3)\gamma^2+(4\mu+3\beta^2\mu-\mu^3)\gamma-\mu^2-\beta^2\mu^2=0.$$ 
Eliminating $\gamma$ from this equation and  $2\gamma^2-\gamma\mu-1=0$, we have
$$\beta^2(2\mu^4+15\mu^2-9)=0,$$
which implies that $\mu$ is a non-zero constant because of $\beta\ne 0$. However, this contradicts \cite[Theorem 1.2]{sa}.
Therefore, $V_{22}$ is an empty set.

{\bf Case (b.2.iii)}:  $f(\gamma,\mu):=
1536\gamma^8+\sum_{i=0}^{7}g_{i}(\mu)\gamma^i=0$ on an open subset $V_{23}\subset V_2$.
 We differentiate  this condition with respect to $e_3$,
 and use (\ref{cd2}), (\ref{eq2.1}) and (\ref{eq2}). Then, putting $\omega=\beta^2$,
  we obtain 
\be
\begin{split}
&\omega(7808\gamma^8-6464\gamma^7\mu-1856\gamma^6\mu^2-1760\gamma^6
+19744\gamma^5\mu^3-2160\gamma^5\mu\\
&-24576\gamma^4\mu^4-2840\gamma^4\mu^2+240\gamma^4+16304\gamma^3\mu^5
+444\gamma^3\mu^3+664\gamma^3\mu\\
&+444\gamma^3\mu^3+664\gamma^3\mu-5826\gamma^2\mu^6-1224\gamma^2\mu^4
+484\gamma^2\mu^2+939\gamma\mu^7\\
&+66\gamma\mu^5+158\gamma\mu^3-51\mu^8+54\mu^6+14\mu^4)\\
&+7808\gamma^{10}-26560\gamma^9\mu+48256\gamma^8\mu^2-5664\gamma^8-59296\gamma^7\mu^3\\
&+12080\gamma^7\mu+50976\gamma^6\mu^4-17256\gamma^6\mu^2+1120\gamma^6
-31888\gamma^5\mu^5\\
&+18356\gamma^5\mu^3+360\gamma^5\mu+13998\gamma^4\mu^6-11596\gamma^4\mu^4
-960\gamma^4\mu^2\\
&-120\gamma^4-3795\gamma^3\mu^7+6138\gamma^3\mu^5+434\gamma^3\mu^3-208\gamma^3\mu
\\
&+528\gamma^2\mu^8-2511\gamma^2\mu^6-1346\gamma^2\mu^4+90\gamma^2\mu^2-27\gamma\mu^9\\
&+480\gamma\mu^7
+386\gamma\mu^5+200\gamma\mu^3-27\mu^8+6\mu^6+38\mu^4=0.
\end{split}\label{eq9}
\ee
 Computing the resultant of the left hand sides of  (\ref{eq7}) and (\ref{eq9}) with respect to $\omega$,
 we obtain
 \be(2\gamma^2-\gamma\mu-1)\sum_{i=0}^{18}h_{i}(\mu)\gamma^i=0,\nonumber\ee
 where $h_i(\mu)$ are
   polynomials given by
  \be
  \begin{split}
   h_0=&1377\mu^{19}-4527\mu^{17}+1284\mu^{15}+4728\mu^{13}\\
   &+1468\mu^{11}-4908\mu^9,\\
  	 h_1=&
  	-66060\mu^{18}+143388\mu^{16}+45504\mu^{14}-74376\mu^{12}\\
  	&-90176\mu^{10}-24512\mu^8,\\
  	h_2=&1397709\mu^{17}-1823607\mu^{15}-1459380\mu^{13}-146496\mu^{11}\\
  	&+289596\mu^9+119828\mu^7,\\
  	h_3=&-17308746\mu^{16}+11826810\mu^{14}+14566216\mu^{12}+6035952\mu^{10}\\
  	&+2043336\mu^8+763816\mu^6,\\
  	 h_4=&140724708\mu^{15}-40031364\mu^{13}-73716152\mu^{11}-29812064\mu^9\\
  	 &-6305728\mu^7+511024\mu^5,\\
  	 h_5=&-801068376\mu^{14}+56606496\mu^{12}+232622128\mu^{10}+61803936\mu^8\\
  	 &-3732032\mu^6-1748096\mu^4,\\
  	h_6=&3336681024\mu^{13}+31432448\mu^{11}-555418672\mu^9-78779056\mu^7\\
  	&+7995776\mu^5-2476096\mu^3,\\
  	h_7=&-10529445888\mu^{12}-118321664\mu^{10}+1089457312\mu^8+68763808\mu^6\\
  	&+4847488\mu^4-430976\mu^2,\\
  	h_8=&25909096832\mu^{11}-541759232\mu^9-1652978624\mu^7+6159040\mu^5\\
  	&+14641152\mu^3+625920\mu,\\
  	h_9=&-50856105728\mu^{10}+3160585216\mu^8+1852903808\mu^6-128791552\mu^4\\
  	&+12595200\mu^2
  	+230400,\\
  	h_{10}=&80930532864\mu^9-8199388160\mu^7-1406354176\mu^5+71124736\mu^3\\
  	&-3322880\mu,\\
  	h_{11}=&-105451162624\mu^8+14029552640\mu^6+551636480\mu^4-60403200\mu^2\\
  	&-3379200,\\
  	h_{12}=&112905166848\mu^7-17520074752\mu^5+317690880\mu^3-40262656\mu,\\
  	 h_{13}=&-98991538176\mu^6+16388972544\mu^4-537110528\mu^2+27381760,\\
  	 h_{14}=&70233264128\mu^5-11463987200\mu^3+441262080\mu,\\
  	h_{15}=&-39343300608\mu^4+5679833088\mu^2-109936640,\\
  	h_{16}=&16633511936\mu^3-1843052544\mu,\\
  	h_{17}=&-4831674368\mu^2+243859456,\\
  	h_{18}=&767557632\mu.
  \end{split}\nonumber
  \ee
  
   Since Case (b.2.ii) does not occur, we have $\sum_{i=0}^{18}h_{i}(\mu)\gamma^i=0$.
The resultant $R_2(\mu)$ of $f(\gamma, \mu)$ and $\sum_{i=0}^{18}h_{i}(\mu)\gamma^i$
with respect to $\gamma$ is given by
$$R_2(\mu)=\mu^{36}k(\mu),$$ where $k(\mu)$ is a polynomial in $\mu$ with constant coefficients of degree $116$.
Since the explicit form of $k(\mu)$ is not important for the argument,
we do not list it.
Thus, we deduce that $\mu$ is constant, that is, the mean curvature is constant.
 According to \cite[Theorem 1.2]{sa}, we conclude that 
${V}_{23}$ is  minimal ruled.


Consequently, $M$ is locally obtained by the construction described in Example \ref{ex1}.
The proof is finished.
\qed

 \end{document}